\documentclass[12pt]{amsart}%
\usepackage{amssymb}
\usepackage{amsmath}
\usepackage{amsfonts}
\usepackage{graphicx}%
\providecommand{\U}[1]{\protect\rule{.1in}{.1in}}
\newtheorem{theorem}{Theorem}
\theoremstyle{plain}
\newtheorem{acknowledgement}{Acknowledgement}

\newtheorem{corollary}{Corollary}

\newtheorem{example}{Example}

\newtheorem{proposition}{Proposition}
\newtheorem{remark}{Remark}

\numberwithin{equation}{section}
\begin{document}
\title[OIRI-algebras]{When Is an OIRI-Algebra a Banach Lattice? Structure Theorems and Counterexamples}
\author{Zied Jbeli}
\address{Department of Information Technology, Higher Institute of Technological
Studies of Kef, University Campus, Boulifa 7100 Kef, Tunisia}
\email{zied.jbeli@fsb.ucar.tn}
\author{Mohamed Ali Toumi}
\address{Laboratoire d'Ing\'{e}nierie Math\'{e}matique LIM, Ecole Polytechnique de
Tunisie, D\'{e}partement de Math\'{e}matiques, Facult\'{e} des Sciences de
Bizerte, Universit\'{e} de Carthage, Carthage, Tunisia}
\email{mohamedali.toumi@fsb.ucar.tn}
\subjclass[2000]{06F25, 13C05, 46B42}
\keywords{OIRI-algebra, relatively uniformly complete $f$-algebra, semiprime, finite element.}

\begin{abstract}
This paper investigates the exact conditions under which a relatively
uniformly complete semiprime $f$-algebra whose order ideals are ring ideals
(an OIRI-algebra) admits a norm turning it into a Banach lattice or an
AM-space. While recent literature claimed that this equivalence holds
automatically in the presence of a weak order unit, we construct a concrete
counterexample within a relatively uniformly complete vector lattice showing
that this assertion fails in general. By analyzing the canonical embedding
$\psi:A\rightarrow Orth(A)$, we prove that the equivalence between being an
OIRI-algebra, admitting an AM-space or Banach lattice norm and satisfying
$Orth(A)=$ $Stab(A)$ $=Z(A$) holds under the assumption that $\psi\left(
A\right)  $ is relatively uniformly closed in $Orth(A)$ (alongside a
characterization via finite elements in the unitization $A^{\sharp}$). We also
deduce immediate consequences for unital $f$-algebras and generalize our main
theorem to non-semiprime $f$-algebras via the quotient modulo the nilpotent
radical $N(A)$.

\end{abstract}
\maketitle

\section{Introduction}

The interplay between algebraic and order-theoretic structures is a central
theme in the theory of $f$-algebras, particularly concerning the conditions
under which order ideals inherently constitute ring ideals. This structural
equivalence was investigated for semiprime $f$-algebras by \textit{Basly} and
\textit{Triki} [2] in 1988, and subsequently expanded to the setting of
general $f$-rings by \textit{Henriksen,} \textit{Larson,} and \textit{Smith}
[5], who introduced the designation OIRI-rings (Order Ideals are Ring Ideals)
for spaces exhibiting this property. More recently, \textit{Boulabiar and
Souissi} [3] advanced the theory by examining OIRI-algebras through the
framework of submultiplicative lattice seminorms. Among other findings, they
established that a semiprime lattice algebra possesses the order-to-ring
property if and only if it admits a submultiplicative lattice norm [3,
Corollary 4.4]. Motivated by this remarkable result, a natural question
arises: \textit{under what exact conditions does an }OIRI\textit{-algebra
admit a lattice norm that makes it a Banach lattice?}

The principal objective of the present paper is to furnish new
characterizations of OIRI-algebras and to resolve structural inaccuracies
present in recent literature. Previous work (specifically [16, Theorems 5])
asserted strict equivalences between OIRI-algebras and Banach lattices under
the assumption that the space is a relatively uniformly complete semiprime
$f$-algebra with a weak order unit. We prove that this assumption is
fundamentally insufficient.

By focusing on the canonical embedding of a relatively uniformly complete
$f$-algebra $A$ into its orthomorphism space $Orth(A),$ we provide a correct
reformulation of these theorems. We establish that the structural equivalence
requires the crucial condition that the canonical image of $A$ is relatively
uniformly closed in $Orth(A).$ To substantiate this framework, we construct a
definitive counterexample inside a universally complete vector lattice,
proving that a relatively uniformly\ complete semiprime OIRI-algebra can
entirely fail to be a Banach lattice when this closure condition is not met.
Furthermore, we extend our main structural theorem beyond the semiprime
setting, showing that the equivalence holds for arbitrary $f$-algebras by
passing to the semiprime quotient algebra $A/N(A),$ where $N(A)$ is the set of
all nilpotent elements in $A.$

\section{Preliminaries}

Recall that a \textit{vector lattice} is a \textit{partially ordered vector
space} in which every pair of elements has a supremum and an infimum. A vector
lattice $A$ is called \textit{Archimedean} if for any $f\in A$ and any $0\leq
g\in A,$ $nf\leq g$ for all $n\in%
\mathbb{N}
$ implies $f\leq0$. Throughout this paper, all vector lattices under
consideration are assumed to be Archimedean and $A^{+}=\left\{  f\in
A:f\geq0\right\}  $ denotes the positive cone of $A.$

Recall that a norm $\left\Vert \cdot\right\Vert $ on a vector lattice A is a
lattice norm if $\left\vert f\right\vert \leq\left\vert g\right\vert $ implies
$\left\Vert f\right\Vert \leq\left\Vert g\right\Vert $ for all $f,g\in A.$ A
vector lattice complete with respect to a lattice norm is a \textit{Banach
lattice}. A Banach lattice $A$ is called an \textit{AM-space} if its norm is
an M-norm, meaning $\left\Vert f\vee g\right\Vert =\max\left\{  \left\Vert
f\right\Vert ,\left\Vert g\right\Vert \right\}  $ for all $f,g\in A^{+}$. A
sublattice $I$ of a vector lattice $A$ is termed an \textit{order ideal} if it
follows from $\left\vert g\right\vert \leq\left\vert f\right\vert $ and $f\in
I$ that $g\in I$. An order ideal $I$ is \textit{order-dense} if and only if
its disjoint complement is trivial (i.e., $I^{d}=\{0\}$). The \textit{order
ideal} generated by a single element $f$ is denoted by $A_{f}.$ An order ideal
$I$ is called a \textit{band} if every subset $J$ of $I$ such that
$\sup(J)=f\in A$ implies that $f\in I.$ An element $e\in A^{+}$ is named a
\textit{weak} \textit{order unit} whenever $g\in A$ and $\left\vert
g\right\vert \wedge e=0$ yield $g=0,$ and it is called a \textit{strong}
\textit{order unit} if for all $g\in A$ there exists $\alpha\in%
\mathbb{R}
^{+}$ such that $|g|\leq\alpha e.$ A set $\{f_{\lambda}:\lambda\in\Lambda\}$
of nonzero pairwise disjoint elements of $A^{+}$ is called a \textit{disjoint
system} in $A$. The disjoint system $\{f_{\lambda}:\lambda\in\Lambda\}$ in $A$
is said to be \textit{maximal} whenever $f\perp f_{\lambda}$ for all
$\lambda\in\Lambda$ implies that $f=0$ (see [10, Definition 28.4]). A vector
lattice\textit{\ }$A$ is named \textit{Dedekind complete} if for each nonempty
order bounded set $B\subset A$, $\sup$ $B$ exists in $A,$ and it is called
\textit{laterally complete} whenever every disjoint system in $A$ has a
supremum in $A$\textit{.} If $A$ is Dedekind complete and laterally complete,
then $A$ is said to be\textit{\ universally complete}. We refer the reader to
the standard monographs [1, 10].

If $A^{+}$ is closed under multiplication, then $A$ is called a
\textit{lattice algebra }(shortly, $\ell$\textit{-algebra}). The set of
\textit{nilpotent elements} of $A$ is denoted by $N(A),$ and is a band of $A.$
Moreover, if $N(A)=\{0\},$ then $A$ is called \textit{semiprime.} The space
$N(A)^{d}$ is the \textit{disjoint complement} of the nilpotent ideal $N(A).$
That is%
\[
N(A)^{d}=\{x\in A:\left\vert x\right\vert \wedge\left\vert y\right\vert
=0\text{ for all }y\in N(A)\}.
\]
An $\ell$-algebra is said to be\textit{ }$f$\textit{-algebra} whenever
$f\wedge g=0$ yields $(fh)\wedge g=(hf)\wedge g=0$ for all $h\in A^{+}.$ If
$A$ is an $f$-algebra, then $A/N(A)$ is an Archimedean semiprime $f$-algebra.
A semiprime $f$-algebra $A$ is said to have the \textit{Stone condition}
whenever $Id_{A}\wedge f\in A$ for all $f\in A$, where $Id_{A}$ is the
identity operator on $A$. For more details, see [4].

To properly contextualize the algebraic structures discussed later, we recall
the notion of orthomorphisms. An order bounded linear operator $\pi
:A\longrightarrow A$ on a vector lattice $A$ is called an
\textit{orthomorphism} whenever $\left\vert f\right\vert \wedge\left\vert
g\right\vert =0$ implies $\left\vert \pi(f)\right\vert \wedge\left\vert
g\right\vert =0$ for all $f,g\in A.$ The collection of all orthomorphisms on
$A$ is denoted by $Orth(A).$ The \textit{center} of $A,$ denoted by $Z(A)$, is
defined as the principal order ideal generated by $Id_{A}$ in $Orth(A).$ A
linear operator $\pi$ on $A$ is \textit{ideal-preserving} if for any order
ideal $I$ of $A$, $\pi(g)\in I$ whenever $g\in I$. We call the
\textit{stabilizer} of $A,$ the set of all ideal-preserving operators on $A$.
This space is denoted by $Stab(A).$\ Furthermore, if $A$ is an $f$-algebra,
every element $f\in A$ induces a canonical multiplication operator $\pi_{f}\in
Orth(A)$, which is defined by $\pi_{f}(g)=fg$ for all $g\in A.$ For more
details, see [4, 18].

For $u\in A^{+},$ a sequence $(f_{n})_{n\in%
\mathbb{N}
}$ in a vector lattice $A$ is said to be $u$\textit{-relatively uniformly
convergent} to $f$ (shortly, $f_{n}\longrightarrow f$ $(u)$) whenever for
every $\epsilon>0,$ there exists $n_{\epsilon}\in%
\mathbb{N}
$ such that $\left\vert f_{n}-f\right\vert \leq\epsilon u$ for all $n\geq
n_{\epsilon}.$ If $f_{n}\longrightarrow f$ $(u)$ for some $u\in A^{+},$ then
$(f_{n})_{n\in%
\mathbb{N}
}$ is called \textit{relatively uniformly convergent} to $f$ (shortly,
$f_{n}\longrightarrow f$ $(ru)$). In addition, relatively uniform limits are
unique, whenever $A$ is Archimedean. A sequence $(f_{n})_{n\in%
\mathbb{N}
}$ is called \textit{relatively uniform Cauchy} if the double sequence
$(f_{n}-f_{m})_{(n,m)\in%
\mathbb{N}
\times%
\mathbb{N}
}$ is relatively uniformly convergent to zero. Next, a vector lattice is said
to be\textit{ relatively uniformly complete} (shortly, $(ru)$%
\textit{-complete}) if every relatively uniform Cauchy sequence is relatively
uniformly convergent. We say that a subset $D$ of $A$ is \textit{relatively
uniformly} \textit{closed} (abbreviated $(ru)$\textit{-closed}) if for each
sequence $\left(  f_{n}\right)  _{n\in%
\mathbb{N}
}$ in $D$ with $f_{n}\longrightarrow f$ $(r.u)$\ we have $f\in D.$ The
\textit{relative uniform topology }is defined by collecting the relatively
uniformly closed sets in $A$. For more details, see [10].

Finite elements in Archimedean vector lattices are introduced by \textit{Weber
}[17] and studied recently by many authors [8, 9, 12, 13, 15, 16]. An element
$f$ in a vector lattice $A$ is named \textit{finite}, if there exists an
element $h\in A^{+}$ such that for every $g\in A$ there is a number
$\lambda_{g}>0$ satisfying the inequality $|g|\wedge n|f|\leq\lambda_{g}h$ for
all $n\in%
\mathbb{N}
$. The element $h$ is named an $A$-majorant of $f$. The set of all finite
elements of $A$ is denoted by $\Phi_{1}(A).$ For more details, see [17].

\section{Main results}

Recall from [5] that an $\ell$-algebra $A$ is called an OIRI-Algebra, whenever
each order ideal of $A$ is a ring ideal. For any element $f\in A,$ let
$\pi_{f}:A\longrightarrow A$ denote the left multiplication operator defined
by $\pi_{f}(x)=fx$ for all $x\in A.$ Throughout $A$ will be an $f$-algebra, we
define the canonical mapping%
\[%
\begin{array}
[c]{cccc}%
\psi: & A & \longrightarrow & Orth(A)\\
& f & \mapsto & \pi_{f}%
\end{array}
\]
which embeds $A$ as a vector sublattice of $Orth(A),$ provided that $A$ is semiprime.

In 1988, \textit{Basly} and \textit{Triki} [2] proved the following result.

\begin{theorem}
[{[2, Theorem 3]}]Let $A$ be an $f$-algebra. Then the following assertions are
equivalent:\newline$\left(  i\right)  $ $A$ is an $OIRI$-algebra,\newline%
$\left(  ii\right)  $ every order ideal in $A$ is a subalgebra of
$A,\newline\left(  iii\right)  $ $\psi(A)$ is a subset of $Z(A).$
\end{theorem}

According to [4, Proposition 12.1], if $A$ is a semiprime $f$-algebra then
$\psi$ is an injective lattice homomorphism. Moreover, a semiprime $f$-algebra
$A$ is lattice and ring isomorphic to a ring ideal and vector sublattice
$\psi(A)$ in $Orth(A)$ (see [4, Page 81]). Throughout, we identify $A$ with
its canonical image $\psi(A).$

Let $A$ be a non-unital semiprime $f$-algebra satisfying the Stone condition,
a property automatically inherited whenever $A$ is $(ru)$-complete (see [7,
Theorem 2.5]). A unit element can be adjoined to $A$ within $Orth(A)$ by
defining%
\[
A^{\sharp}=A\oplus%
\mathbb{R}
Id_{A}=\{g+\alpha Id_{A}:g\in A~\text{and}~\alpha\in%
\mathbb{R}
\},
\]
where $Id_{A}\notin A$ denotes the identity operator on $A.$
\textit{Huijsmans} and \textit{de Pagter} established that $A^{\sharp}$ forms
a unital $f$-subalgebra of $Orth(A)$ containing $A$ (see [6, Section 7]).

We recall the following useful result from the recent paper [16].

\begin{proposition}
[{[16, Theorem 7]}]Let $A$ be an $(ru)$-complete semiprime $f$-algebra and let
$f\in A^{+}.$ Then the listed assertions are equivalent:\newline$\left(
i\right)  $ $\pi_{f}\in\Phi_{1}(A^{\sharp})$,\newline$\left(  ii\right)  $ the
principal order ideal $A_{f}$ generated by $f$ in $A$ is a ring ideal.
\end{proposition}

As an immediate consequence of Proposition 1, we give the following
characterization of OIRI-algebras via finite elements in the unital
$f$-algebra $A^{\sharp}.$

\begin{theorem}
Let $A$ be an $(ru)$-complete semiprime $f$-algebra. Then the following
statements are equivalent:\newline$\left(  i\right)  $ $A$ is an
OIRI-algebra,\newline$\left(  ii\right)  $ $\pi_{f}\in\Phi_{1}(A^{\sharp})$
for all $f\in A^{+}.$
\end{theorem}

Recently, \textit{Boulabiar and Souissi }established that a semiprime $\ell
$-algebra $A$ is an OIRI-algebra if and only if it admits a submultiplicative
lattice norm (see [3, Corollary 4.4]). However, the existence of such a norm
does not guarantee that the algebra can be endowed with a complete topological
structure. In general, an OIRI-algebra might not admit any Banach lattice
norm, even under the assumption that the $\ell$-algebra is $(ru)$-complete.
The following example explicitly illustrates this claim.

\begin{example}
Let $B$ be a universally complete vector lattice admitting a countably
infinite maximal disjoint system $S=\{e_{n}:n=1,2,\ldots\}\subset B^{+}.$ (For
instance, one may take $B=\left(
\mathbb{R}
^{%
\mathbb{N}
}\right)  ^{u},$ the universally complete space of all real sequences). By the
universal completeness of $B,$ the supremum $e=\underset{n=1}{\overset{\infty
}{\bigvee}}e_{n}$ exists in $B.$ Furthermore, $B$ can be uniquely equipped
with an $f$-algebra multiplication under which $e$ serves as the
multiplicative identity element (see [14, Lemma 3.7]). Under this
multiplication, the elements $\{e_{n}\}_{n=1}^{\infty}$ are pairwise disjoint
idempotents, that is, $e_{n}^{2}=e_{n}$ and $e_{j}e_{k}=0$ for all $j\neq k.$
For each $n\in%
\mathbb{N}
,$ define $f_{n}=\underset{k=1}{\overset{n}{\sum}}e_{k},$ and consider the
vector space:%
\[
A=\underset{n=1}{\overset{\infty}{\bigcup}}B_{f_{n}}%
\]
where $B_{f_{n}}$ is the principal order ideal generated by $f_{n}$ in $B.$
Since each $B_{f_{n}}$ is an order ideal of $B,$ $A$ is an order ideal of the
universally complete space $B,$ and is therefore $(ru)$-complete. To verify
that $A$ is closed under multiplication, let $x,y\in A^{+}.$ Then there exist
$n_{0}\in%
\mathbb{N}
$ and $\lambda_{1},\lambda_{2}>0$ such that $0\leq x\leq\lambda_{1}f_{n_{0}}$
and $0\leq y\leq\lambda_{2}f_{n_{0}}.$ Since the elements $e_{k}$ are pairwise
disjoint idempotents, we have $f_{n_{0}}^{2}=f_{n_{0}}$, which yields%
\[
0\leq x\cdot y\leq\left(  \lambda_{1}\lambda_{2}\right)  f_{n_{0}}^{2}=\left(
\lambda_{1}\lambda_{2}\right)  f_{n_{0}}%
\]
Thus $(A,\cdot)$ is an $(ru)$-complete semiprime $f$-algebra. Furthermore, for
any $x\in A,$ there exists $n_{0}\in%
\mathbb{N}
$ and $M>0$ such that%
\[
\left\vert x\right\vert \leq Mf_{n_{0}}\leq Me.
\]
For any $y\in A,$ we have
\[
\left\vert \pi_{x}(y)\right\vert =\left\vert x\cdot y\right\vert \leq\left(
Me\right)  \cdot\left\vert y\right\vert .
\]
Accordingly, $\pi_{x}\in Z(A).$ Hence $\Psi(A)\subset Z(A).$ In view of
Theorem 1, $A$ is an OIRI-algebra. We now claim that $A$ does not admit any
lattice norm under which it is a Banach lattice. Suppose, per contra, that
there exists a lattice norm $\left\Vert \cdot\right\Vert $ under which $A$ is
a Banach lattice. Consider the infinite series:%
\[
u=\underset{k=1}{\overset{\infty}{\sum}}\frac{2^{-k}}{1+\left\Vert
f_{k}\right\Vert }f_{k}%
\]
Since $\left\Vert \frac{2^{-k}}{1+\left\Vert f_{k}\right\Vert }f_{k}%
\right\Vert \leq2^{-k},$ the series is absolutely convergent in $(A,\left\Vert
\cdot\right\Vert ).$ Moreover, as $A$ is complete, we have $u\in A.$ Hence
there exists $\lambda>0$ and $n_{0}\in%
\mathbb{N}
$ such that $u\leq\lambda f_{n_{0}}.$ On the other hand, since all the terms
in the series are strictly positive, it follows that
\[
u\geq\frac{2^{-(n_{0}+1)}}{1+\left\Vert f_{n_{0}+1}\right\Vert }f_{n_{0}+1}.
\]
Hence%
\[
\frac{2^{-(n_{0}+1)}}{1+\left\Vert f_{n_{0}+1}\right\Vert }f_{n_{0}+1}%
\leq\lambda f_{n_{0}}.
\]
Consequently
\[
f_{n_{0}+1}\leq2^{n_{0}+1}\left(  1+\left\Vert f_{n_{0}+1}\right\Vert \right)
\lambda f_{n_{0}}.
\]
Let $C=2^{n_{0}+1}\left(  1+\left\Vert f_{n_{0}+1}\right\Vert \right)
\lambda,$ we have $f_{n_{0}}+e_{n_{0}+1}\leq Cf_{n_{0}}.$ Let $P_{e_{n_{0}+1}%
}$ denote the band projection onto the principal band generated by
$e_{n_{0}+1}$ in $B.$ Applying $P_{e_{n_{0}+1}}$ to both sides of the
inequality and observing that $P_{e_{n_{0}+1}}(f_{n_{0}+1})=0,$ we obtain
$e_{n_{0}+1}\leq0,$ which contradicts the fact that $e_{n_{0}+1}>0.$
\end{example}

To characterize the coincidence of the multiplier algebras
$Orth(A)=Z(A)=Stab(A),$ it was claimed in [16, Theorem 5] if $A$ is an
$(ru)$-complete semiprime $f$-algebra with a weak order unit, then $A$ is a
Banach lattice if and only if $\pi_{f}\in\Phi_{1}(A^{\sharp})$ for all $f\in
A^{+}.$ However, this assertion is mathematically incorrect. Because their
subsequent main result ([16, Theorem 6]) relies intrinsically on this
equivalence, its general validity is consequently compromised. We present the
following concrete counterexample to explicitly prove why this claimed
equivalence fails.

\begin{example}
Let $A=C_{0}\left(  \mathbb{R}\right)  $ be the vector lattice of all
continuous functions vanishing at infinity with the supremum norm and
pointwise multiplication. It's not hard to show that $A$ is a Banach semiprime
$f$-algebra. Let $f_{n}$ in $A$ be defined by $f_{n}\left(  x\right)
=e^{-\frac{x^{2}}{n+1}},$ for all $n\in\mathbb{N}.$ It follows that $f_{n}\leq
f_{n+1},$ for all $n\in\mathbb{N}.$ and $f_{0}$ is a weak order unit of $A.$
Consider the vector space:%
\[
B=\underset{n=0}{\overset{\infty}{\bigcup}}A_{f_{n}}%
\]
where $A_{f_{n}}$ is the principal order ideal generated by $f_{n}$ in $A.$
Since each $A_{f_{n}}$ is an order ideal of $A,$ $B$ is an order ideal of $A,$
and is therefore $(ru)$-complete. Moreover, $B$ is a semiprime $f$-algebra and
$B$ is an $OIRI$-algebra. By applying Theorem 2, $\pi_{g}\in\Phi_{1}%
(B^{\sharp})$ for all $g\in B^{+}.$ We now claim that $A$ does not admit any
lattice norm under which it is a Banach lattice. Suppose, per contra, that
there exists a lattice norm $\left\Vert \cdot\right\Vert $ under which $A$ is
a Banach lattice. Consider the infinite series:%
\[
u=\underset{n=0}{\overset{\infty}{\sum}}\frac{2^{-n}}{1+\left\Vert
f_{n}\right\Vert }f_{n}%
\]
Since $\frac{2^{-n}}{1+\left\Vert f_{n}\right\Vert }\leq2^{-n},$ the series is
absolutely convergent in $(B,\left\Vert \cdot\right\Vert ).$ Moreover, as $B$
is complete, we have $u\in B.$ Hence there exists $\lambda>0$ and $n_{0}\in%
\mathbb{N}
$ such that $u\leq\lambda f_{n_{0}}.$ Consequently
\[
0\leq\frac{2^{-n}}{1+\left\Vert f_{n}\right\Vert }f_{n}\leq\lambda f_{n_{0}}%
\]
for all $n\in\mathbb{N}.$ By choosing $n$ such that $n>n_{0},$ we deduce that%
\[
0\leq\frac{2^{-n}e^{-\frac{x^{2}}{n+1}}}{1+\left\Vert f_{n}\right\Vert }%
\leq\lambda e^{-\frac{x^{2}}{n_{0}+1}}%
\]
for all $x\in\mathbb{R}.$ Therefore
\[
0\leq\frac{2^{-n}e^{x^{2}(\frac{1}{n_{0}+1}-\frac{1}{n+1})}}{1+\left\Vert
f_{n}\right\Vert }\leq\lambda
\]
As $x\rightarrow+\infty,$ we get
\[
0\leq+\infty\leq\lambda
\]
a contradiction.
\end{example}

Now we prove the first main result of our paper.

\begin{theorem}
Let $A$ be an $(ru)$-complete semiprime $f$-algebra such that its canonical
image $\psi(A)$ is $(ru)$-closed in $Orth(A)$. Then the statements that follow
are equivalent:\newline$\left(  i\right)  $ $A$ is an OIRI-algebra,\newline%
$\left(  ii\right)  $ $\pi_{f}\in\Phi_{1}(A^{\sharp})$ for all $f\in A^{+},$
\newline$\left(  iii\right)  $ There exists an $M$-norm under which $A$ is an
$AM$-space,\newline$\left(  iv\right)  $ There exists a lattice norm under
which $A$ is a Banach lattice,\newline$\left(  v\right)  $
$Orth(A)=Stab(A)=Z(A).$
\end{theorem}

\begin{proof}
$(i)\Leftrightarrow(ii)$ Follows from Theorem 2.

$(i)\Rightarrow(iii)$ Assume that $A$ is an OIRI-algebra. It follows from
Theorem 1 that $\psi(A)\subset Z(A).$ Since $A$ is $(ru)$-complete, it follows
from [4, Theorem 9.5] that $Z(A)$ is also $(ru)$-complete. Moreover, as
$Z(A)=Orth(A)_{Id_{A}},$ then $Z(A)$ equipped with the gauge function of
$[-Id_{A},Id_{A}]$ given by
\[
\rho_{Id_{A}}(\pi)=\inf\{\lambda\in%
\mathbb{R}
:-\lambda Id_{A}\leq\pi\leq\lambda Id_{A}\}
\]
which is an $M$-norm, is an $AM$-space (see [11, Chap II, Proposition 7.2]).
By assumption, $\psi(A)$ is $(ru)$-closed in $Orth(A)$. It is not hard to
verify that $\psi(A)$ is also $(ru)$-closed in $Z(A)$. Recall that norm
convergence and $(ru)$-convergence coincide in any $AM$-space with a strong
order unit, $\psi(A)$ is equivalently norm-closed. Thus $\psi(A)$ is a
complete subspace under the $M$-norm $\rho_{Id_{A}}$. Restricting this norm to
$A$ establishes that $A$ is an $AM$-space.

$(iii)\Rightarrow(iv)$ Obvious.

$(iv)\Rightarrow(v)$ Follows from [18, Corollary 4.2].

$(v)\Rightarrow(i)$ Since $Orth(A)=Stab(A)=Z(A),$ then the canonical embedding
satisfies $\psi(A)\subset Orth(A)=Z(A).$ By Theorem 1, this inclusion
guarantees that $A$ is an OIRI-algebra, and we are done.
\end{proof}

\begin{remark}
1) The assumption $"\psi(A)$ is $(ru)$-closed in $Orth(A)"$ is only used to
prove $(i)\Rightarrow(iii),$ and is not superfluous in Theorem 3. To see this,
we consider the vector lattice $A$ in the above Example 1. We already showed
that $A$ is an OIRI-algebra but does not possess a lattice norm that makes it
as a Banach lattice. We claim that the canonical image $\psi(A)=\left\{
\pi_{x}:x\in A\right\}  $ is not $(ru)$-closed in $Orth(A).$ Indeed, since the
vector lattice $B$ is universally complete, the disjoint sequence
$\{ke_{k}\}_{k=1}^{\infty}$ admits a supremum $f=\underset{k\in%
\mathbb{N}
}{\bigvee}ke_{k}\in B.$ Recall that $e=\underset{k=1}{\overset{\infty}%
{\bigvee}}e_{k}$ and $f_{n}=\underset{k=1}{\overset{n}{\bigvee}}e_{k}$. For
any $n\in%
\mathbb{N}
,$ we have
\[
e-f_{n}=\underset{k=n+1}{\overset{\infty}{\bigvee}}e_{k}\leq\frac{1}%
{n}\underset{k=n+1}{\overset{\infty}{\bigvee}}ke_{k}\leq\frac{1}{n}f.
\]
For any $x\in A^{+},$
\[
0\leq e\cdot x-f_{n}\cdot x\leq\frac{1}{n}f\cdot x
\]
It follows that
\[
0\leq Id_{A}(x)-\pi_{f_{n}}(x)\leq\frac{1}{n}\pi_{f}(x)\text{ for all }x\in
A^{+}.
\]
Next we verify that $\left(  \pi_{f}\right)  _{/A}\in Orth(A).$ For any $x\in
A^{+},$ there exists $n_{0}\in%
\mathbb{N}
$ and $\lambda>0$ such that $x\leq\lambda f_{n_{0}}.$ Inasmuch as the $e_{k}%
$'s are mutually disjoint idempotents,%
\[
f\cdot f_{n_{0}}=\underset{k=1}{\overset{n_{0}}{\bigvee}}ke_{k}\leq
n_{0}f_{n_{0}}%
\]
Consequently,
\[
0\leq f\cdot x\leq\lambda(f\cdot f_{n_{0}})\leq\lambda n_{0}f_{n_{0}}\in A.
\]
Therefore $\pi_{f}(x)\in A$ for all $x\in A,$ and so $\left(  \pi_{f}\right)
_{/A}$ is a well-defined positive orthomorphism in $Orth(A).$ Finally, we
obtain
\[
0\leq Id_{A}-\pi_{f_{n}}\leq\frac{1}{n}\pi_{f}\text{ in }Orth(A).
\]
Thus $\left(  \pi_{f_{n}}\right)  _{/A}\longrightarrow Id_{A}$ $(ru)$ in
$Orth(A).$ Suppose, for a contradiction, that $Id_{A}\in\psi(A).$ Then there
exists $a\in A$ such that $\left(  \pi_{a}\right)  _{/A}=Id_{A}.$ Since $A$ is
an order-dense ideal in $B$ then $\pi_{a}=Id_{B}$ . As $Id_{B}=\pi_{a},$ it
follows that $a=e,$ which implies $e\in A.$ This contradicts the fact that
$e\notin A.$ Therefore $Id_{A}\notin\psi(A).$

2) Theorem 3 provides corrections to both Theorem 5 and Theorem 6 in the paper
[16].\newline3) Uniqueness up to equivalence and Banach f-algebra structure:
Theorem 3 guarantees the existence of an M-norm and a Banach lattice norm on
$A$ whenever $\Psi(A)$ is (ru)-closed in $Orth(A)$. By [1, Theorem 4.3], all
lattice norms under which $A$ is a Banach lattice are mutually equivalent, and
each is equivalent to the restricted gauge norm $\rho_{Id_{A}}(\pi_{f})$.
Moreover, because $\rho_{Id_{A}}$ is submultiplicative on $Z(A)$, $A$
automatically becomes a Banach $f$-algebra (with jointly continuous
multiplication) whenever conditions (i)--(v) hold.
\end{remark}

According to [4, Proposition 12.1], if an $f$-algebra $A$ possesses a
multiplicative identity, its canonical image satisfies $\psi(A)=Orth(A)$.
Consequently, $\psi(A)$ is trivially $(ru)$-closed in $Orth(A).$ Furthermore,
the unitization algebra $A^{\sharp}$ coincides strictly with $A$ itself (i.e.,
$A^{\sharp}=A$). As the global structural hypotheses of Theorem 3 are
intrinsically satisfied under the presence of a multiplicative identity, we
immediately deduce the following corollary.

\begin{corollary}
Let $A$ be an $(ru)$-complete $f$-algebra with a multiplicative identity. Then
the following statements are equivalent:\newline$\left(  i\right)  $ $A$ is an
OIRI-algebra,\newline$\left(  ii\right)  $ $\pi_{f}\in\Phi_{1}(A)$ for all
$f\in A^{+},$ \newline$\left(  iii\right)  $ There exists an $M$-norm under
which $A$ is an $AM$-space,\newline$\left(  iv\right)  $ There exists a
lattice norm under which $A$ is a Banach lattice,\newline$\left(  v\right)  $
$Orth(A)=Stab(A)=Z(A).$
\end{corollary}

\begin{example}
Let $X$ be a completely regular Hausdorff space, and let $C(X)$ be the
$f$-algebra of all continuous real-valued functions on $X.$ We can precisely
characterize when $C(X)$ is an OIRI-algebra based on the topology of $X:$
\[
C(X)\text{ is an OIRI-algebra if and only if }X\text{ is pseudocompact.}%
\]
Assume $X$ is pseudocompact. It follows that every function $f\in C(X)$ is
bounded. Therefore, for any $f\in C(X),$ there exists a real number
$\lambda>0$ such that $\left\vert f\right\vert \leq\lambda1,$ where $1$ is the
constant function one. Hence the multiplicative identity $1$ is a strong order
unit for $C(X).$ Consequently, $C(X)$ equipped with the supremum norm is an
AM-space (and thus a Banach lattice), which guarantees it is an OIRI-algebra
(see [5, Proposition 2.10]). A canonical example is $X=[0,1].$ Conversely,
assume $X$ is not pseudocompact. Then there exists an unbounded function $f\in
C(X)^{+}$ (for instance, $f(x)=x$ in $C(%
\mathbb{R}
)$). Inasmuch as $f$ is unbounded, there is no $\lambda>0$ such that
$f\leq\lambda1.$ Thus the multiplicative identity $1$ fails to be a strong
order unit. If $C(X)$ were a Banach lattice, its multiplicative unit $1$ would
necessarily be a strong order unit. Because it is not, $C(X)$ cannot admit a
Banach lattice norm. In view of Corollary 1, this proves that for any
non-pseudocompact space $X,$ $C(X)$ is not an OIRI-algebra.
\end{example}

Next, motivated by the characterizations of OIRI-rings proved by
\textit{Henriksen et al.} (see [5, Theorem 2.13 and Corollary 2.14]), we give
the following result.

\begin{proposition}
Let $A$ be an $f$-algebra. Then the following statements are
equivalent:\newline$\left(  i\right)  $ $A$ is an OIRI-algebra,\newline%
$\left(  ii\right)  $ The semiprime quotient $A/N(A)$ is an
OIRI-algebra,\newline$\left(  iii\right)  $ $A/N(A)$ can be embedded as an
$f$-subalgebra into $Z(A/N(A)).$\textbf{\newline}
\end{proposition}

\begin{proof}
$(i)\Rightarrow(ii)$\textbf{ }Let $[a],[x]\in(A/N(A))^{+}.$ We can choose
representatives $a,x\in A^{+}.$ Since $A$ is an OIRI-algebra, the principal
order ideal $A_{x}$ is a ring ideal, and hence $ax\leq\lambda x$ for some
$\lambda>0.$ Passing to the quotient gives $[a][x]=[ax]\leq\lambda\lbrack x].$
Thus $A/N(A)$ is an OIRI-algebra.

$(ii)\Rightarrow(i)$ Let $a,x\in A^{+}.$ Since $A/N(A)$ is an OIRI-algebra,
there exists $\lambda>0$ such that $[a][x]\leq\lambda\lbrack x].$ Hence
$[ax-\lambda x]\leq0.$ Consequently, $(ax-\lambda x)^{+}\in N(A).$ Inasmuch as
$ax\in N(A)^{d}$ and $N(A)^{d}$ is an order ideal, the inequality
$0\leq(ax-\lambda x)^{+}\leq ax$ implies $(ax-\lambda x)^{+}\in N(A)^{d}.$
Therefore $(ax-\lambda x)^{+}\in N(A)\cap N(A)^{d}=\{0\}.$ Accordingly,
$ax\leq\lambda x.$ Thus $A$ is an OIRI-algebra.

$(ii)\Rightarrow(iii)$ Let $B=A/N(A),$ which is a semiprime $f$-algebra. For
any $[a]\in B,$ consider the multiplication operator $\pi_{\lbrack
a]}([x])=[a][x].$ Since $B$ is an OIRI-algebra, then $\pi_{\lbrack a]}\in
Z(B).$ The canonical mapping $\varphi:B\longrightarrow Z(B)$ defined by
$\varphi([a])=\pi_{\lbrack a]}$ is an injective ring and lattice homomorphism,
and hence embeds $A/N(A)$ as an $f$-subalgebra into $Z(A/N(A)).$

$(iii)\Rightarrow(ii)$ Suppose that $A/N(A)$ embeds into $Z(A/N(A))$ via
$[a]\longrightarrow\pi_{\lbrack a]}$. It follows that there exists $\lambda>0$
such that $\left\vert \pi_{\lbrack a]}\right\vert \leq\lambda Id_{A/N(A)}$ for
any $[a]\in A/N(A).$ Consequently,
\[
\left\vert \lbrack a][x]\right\vert =\left\vert \pi_{\lbrack a]}\left(
[x]\right)  \right\vert \leq\lambda Id_{A/N(A)}\left(  [x]\right)
=\lambda\lbrack x]
\]
for all $[x]\in(A/N(A))^{+}.$ Therefore the principal order ideal generated by
$[x]$ is a ring ideal. Thus $A/N(A)$ is an OIRI-algebra.
\end{proof}

Recall that a seminorm $p$ on a vector lattice $A$ is called a lattice
seminorm (or Riesz seminorm) if $\left\vert f\right\vert \leq\left\vert
g\right\vert $ implies $p(f)\leq p(g)$ for all $f,g\in A$. A lattice seminorm
$p$ is called M-seminorm if
\[
p(f\vee g)=\max\left\{  p(f),p(g)\right\}  \text{ for all }f,g\in A^{+}.
\]
Furthermore, in the context of an $\ell$-algebra $A$\textbf{, }after Boulabiar
and Souissi [3]\textbf{,} a seminorm $p$ is termed \textit{nil-faithful} if it
vanishes exclusively on nilpotent elements. That is $p(f)=0$ implies $f\in
N(A)$.

Now we prove the second main result of our work.

\begin{theorem}
Let $A$ be an $(ru)$-complete $f$-algebra. Let $B=A/N(A)$ denote its canonical
semiprime quotient, let $q:A\longrightarrow B$ be the canonical quotient map,
and let $\varphi:B\longrightarrow Orth(B)$ be the natural embedding defined by
$\varphi\left(  \left[  f\right]  \right)  =\pi_{\left[  f\right]  },$ where
$\pi_{\left[  f\right]  }\left(  \left[  x\right]  \right)  =\left[  f\right]
\left[  x\right]  .$ Suppose that the canonical image $\varphi(B)$ is
$(ru)$-closed in $Orth(B)$. Then the following statements are
equivalent:\newline$\left(  i\right)  $ $A$ is an OIRI-algebra,\newline%
$\left(  ii\right)  $ $B$ is an OIRI-algebra,\newline$\left(  iii\right)  $
$B$ can be embedded as an $f$-subalgebra in $Z(B),$\newline$\left(  iv\right)
$ There exists an $M$-norm under which $B$ is an AM-space,\newline$\left(
v\right)  $ There exists a lattice norm under which $B$ is a Banach
lattice,\newline$\left(  vi\right)  $ $Orth(B)=Stab(B)=Z(B),\newline$If, in
addition, $\varphi(q(N(A)^{d}))$ is $(ru)$-closed in $Orth(B),$ then the
following statements are also equivalent to the above:$\newline\left(
vii\right)  $ There exists a nil-faithful $M$-seminorm on $A$ such that
$N(A)^{d},$ equipped with the restriction of this $M$-seminorm, is an
AM-space,\newline$\left(  viii\right)  $ There exists a nil-faithful lattice
seminorm on $A$ such that $N(A)^{d},$ equipped with the restriction of this
lattice seminorm, is a Banach lattice.\newline
\end{theorem}

\begin{proof}
First, the equivalence of $(i),$ $(ii)$ and $(iii)$ follows from Proposition 2.

Next, since $A$ is an $(ru)$-complete $f$-algebra, then $B$ is an
$(ru)$-complete semiprime $f$-algebra. Consequently, the equivalence of
$(ii),$ $(iv),$ $(v),$ and $(vi)$ is a consequence of Theorem 3.

Moreover, if $\varphi(q(N(A)^{d}))$ is $(ru)$-closed in $Orth(B),$ we have:

$(iv)\Rightarrow(vii)$ Assume that $B$ is an AM-space under a lattice norm
$\left\Vert \cdot\right\Vert _{B}$. Define $p:A\longrightarrow%
\mathbb{R}
^{+}$ by $p(f)=\left\Vert \left[  f\right]  \right\Vert _{B}$. Since
\[
p(f)=0\text{ if and only if }f\in N(A),
\]
$p$ is a nil-faithful lattice seminorm on $A.$ Restricted to $N(A)^{d},$ $p$
acts as a lattice M-norm. Let $(f_{n})$ be a $p$-Cauchy sequence in $N(A)^{d}%
$. Then $([f_{n}])$ is a Cauchy sequence in $q\left(  N(A)^{d}\right)  $.
Since $B$ is a Banach lattice, $([f_{n}])$ converges in norm to some limit
$\left[  g\right]  \in B.$ Inasmuch as multiplication in the Banach lattice
$B$ is continuous, the canonical embedding $\varphi:B\longrightarrow Orth(B)$
is a bounded linear operator. Therefore $\varphi([f_{n}])\longrightarrow
\varphi(\left[  g\right]  )$ in $Orth(B)$. Since $B$ is a Banach lattice,
$Orth(B)=Z(B)$ is an AM-space with a strong order unit, and so its
$(ru)$-topology and norm topology in $Orth(B)$ coincide. Since $\varphi
(q(N(A)^{d}))$ is $(ru)$-closed in $Orth(B),$ it follows that $\varphi(\left[
g\right]  )\in\varphi(q(N(A)^{d}))$. Consequently, $\varphi(\left[  g\right]
)=\varphi([h])$ for some $h\in N(A)^{d}$. Since $B$ is semiprime, $\varphi$ is
injective, implying $\left[  g\right]  =[h].$ Hence $\left\Vert \left[
f_{n}\right]  -\left[  h\right]  \right\Vert _{B}\longrightarrow0,$ yielding
\[
p(f_{n}-h)=\left\Vert \left[  f_{n}\right]  -\left[  h\right]  \right\Vert
_{B}\longrightarrow0.
\]
Therefore $q(N(A)^{d})$ contains its limit and is closed in $B,$ and hence
$N(A)^{d}$ is complete under $p.$

$\left(  vii\right)  \Rightarrow(viii)$ Obvious.

$(viii)\Rightarrow(ii)$ Assume $N(A)^{d}$ is a Banach space under a
nil-faithful lattice seminorm. The restriction
\[
q|_{N(A)^{d}}:N(A)^{d}\longrightarrow q(N(A)^{d})
\]
is a lattice isomorphism, making $q(N(A)^{d})$ a Banach lattice under the
induced norm $\left\Vert \left[  f\right]  \right\Vert _{B}=p(f)$ for all
$f\in N(A)^{d}$
\[
Orth(q(N(A)^{d}))=Z(q(N(A)^{d})).
\]
Let $\left[  a\right]  ,\left[  b\right]  \in B_{+}.$ As $ab,a^{2},b^{2}\in$
$N(A)^{d},$ we deduce that there exists $\lambda>0$ such that
\[
0\leq\left(  \left[  a\right]  \left[  b\right]  \right)  ^{2}=\pi_{\left[
b\right]  ^{2}}\left(  \left[  a\right]  ^{2}\right)  \leq\lambda\left[
a\right]  ^{2}.
\]
Since $B$ is a semiprime $f$-algebra, it follows that
\[
0\leq\left[  a\right]  \left[  b\right]  \leq\sqrt{\lambda}\left[  a\right]  .
\]

Therefore the order ideal generated by $\left[  a\right]  $ in $B$ is a ring
ideal and so $B$ is an OIRI-algebra.
\end{proof}

\begin{remark}
This theorem serves as a generalization of previous characterizations of
OIRI-algebras and Banach lattices, seamlessly recovering the known results in
the more restrictive cases where $A$ is semiprime or unital.

1. The Semiprime Case: Suppose $A$ is a semiprime $f$-algebra. By definition,
the nilpotent ideal is trivial (i.e., $N(A)=\{0\}$). In this case, the
quotient space $B=A/N(A)$ is canonically isomorphic to $A$ itself, and the
quotient map $q:A\longrightarrow B$ acts as the identity map. Consequently,
the disjoint complement is $N(A)^{d}=\{0\}^{d}=A.$ Since $q(N(A)^{d})=B,$ the
additional hypothesis that $\varphi(q(N(A)^{d}))$ is $(ru)$-closed becomes
completely redundant, as it is identical to the primary hypothesis on
$\varphi(B).$ Furthermore, since $\ker(p)=N(A)=\{0\},$ the nil-faithful
seminorms in statements $(vii)$ and $(viii)$ tighten into strict lattice norms
on all of $A.$ Thus the theorem collapses exactly into Theorem 3.

2. The Unital Case: Suppose $A$ is an $f$-algebra possessing a multiplicative
identity $e.$ It is a known fact that any unital $f$-algebra is automatically
semiprime, so again $N(A)=\{0\}$ and $B=A.$ More specifically, the presence of
a multiplicative identity guarantees that the natural embedding $\varphi
:A\longrightarrow Orth(A)$ is a lattice isomorphism (i.e., $Orth(A)\cong A$).
As the image $\varphi(A)$ is the entirety of $Orth(A),$ it is trivially
$(ru)$-closed. Thus both hypotheses of the theorem are automatically
satisfied. Therefore the theorem perfectly mirrors Corollary 1.
\end{remark}

\begin{remark}
Sharpness of the restriction to $N(A)^{d}$: In statements (vii) and (viii) of
Theorem 4, completeness is established specifically on the disjoint complement
$N(A)^{d}$ rather than on the whole algebra $A$. This restriction is essential
and cannot be strengthened to conclude that $A$ itself admits a Banach lattice
norm. To see this, consider the Archimedean vector lattice $A=\mathbb{R}%
\oplus\mathbb{R}^{\mathbb{N}}$ endowed with the coordinatewise order and the
multiplication defined by%
\[
(\alpha,x)\text{\textperiodcentered}(\beta,y)=(\alpha\beta,0)\text{ for all
}\alpha,\beta\in\mathbb{R}\text{ and }x,y\in\mathbb{R}^{\mathbb{N}}.
\]

Clearly, $A$ is an $(ru)$-complete $f$-algebra (as a direct sum of two
(ru)-complete vector lattices) whose nilpotent radical is the band
$N(A)=\{0\}\oplus\mathbb{R}^{\mathbb{N}}$ and whose disjoint complement is
$N(A)^{d}=\mathbb{R}\oplus\{0\}$. Moreover, every order ideal $I$ of $A$
decomposes as $I=I_{1}\oplus I_{2}$, where $I_{1}$ and $I_{2}$ are order
ideals of $\mathbb{R}$ and $\mathbb{R}^{\mathbb{N}}$, respectively; under the
defined product, $I$ is manifestly a ring ideal, so $A$ is an OIRI-algebra.
The semiprime quotient $B=A/N(A)$ is canonically isomorphic to $\mathbb{R}$,
which is a unital AM-space. Consequently, $\varphi(B)=Orth(B)$ and
$\varphi(q(N(A)^{d}))=\varphi(B)$ are trivially (r.u)-closed in $Orth(B)$, so
all hypotheses of Theorem 4 are satisfied. Defining $p(\alpha,x)=|\alpha|$
yields a nil-faithful M-seminorm on $A$ under which $N(A)^{d}=\mathbb{R}%
\oplus\{0\}$ is an AM-space, confirming (vii) and (viii). However, $A$ itself
cannot admit any lattice norm turning it into a Banach lattice, because its
closed order ideal $\{0\}\oplus\mathbb{R}^{\mathbb{N}}$ is an
infinite-dimensional universally complete vector lattice and therefore admits
no complete lattice norm (as shown by the series argument in Example 1).
\end{remark}

\begin{acknowledgement}
The first author would like to express his sincere gratitude to the anonymous
referee for his careful reading of an earlier preliminary version of this
manuscript. Through an illustrative example in his report, the referee showed
that, for an (ru)-complete semiprime f-algebra with a weak order unit, the
equivalence between the OIRI-algebra property and the Banach lattice structure
does not hold in full generality. This insightful observation directly
motivated the refined theoretical framework---as well as the explicit
counterexamples constructed in Section 3---established in the present work.
\end{acknowledgement}

\section{Declarations}

\textbf{Author Contributions}

Both authors contributed equally to this work.

\textbf{Funding}

None.

\textbf{Data availability}

No datasets were generated or analyzed during the current study.

\textbf{Ethical Approval}

Not applicable.

\textbf{Competing Interests}

The authors declare no competing interests.


\begin{thebibliography}{99}                                                                                               %


\bibitem {1}Aliprantis, C.D.: Burkinshaw, O.: Positive Operators. Academic
Press Inc., London (1985).

\bibitem {2}Basly, M., Triki, A.: $F$-algebras in which order ideals are ring
ideals, Proc. Konin. Neder. Akad. Wet. 91, 231--234 (1988).

\bibitem {3}Boulabiar, K., Souissi, R.: Order-to-ring ideals via
submultiplicative lattice seminorms in lattice-ordered algebras,
arXiv:2606.20618v1 (2026).

\bibitem {4}De Pagter, B.: $f$-algebras and Orthomorphisms, PhD. Thesis,
Leiden (1981).

\bibitem {5}Henriksen, M., Larson, S., and Smith, F. A.: When is every order
ideal a ring ideal?. Comment. Math. Univ. Carolin. 32, (3), 411-416 (1991).

\bibitem {6}Huijsmans, C.B., de Pagter, B.: Averaging operators and positive
contractive projections. J. Math. Anal. Appl. 113, 163--184 (1986).

\bibitem {7}Huijsmans, C.B., de Pagter, B.: The order bidual of
lattice-ordered algebras. J. Funct. Analysis 59, 41--64 (1984).

\bibitem {8}Jbeli, Z., Toumi, M.A.: Characterizations of self-majorizing
elements in Archimedean vector lattices. Algebra Univers. 85, 29 (2024).

\bibitem {9}Jbeli, Z., Toumi, M.A. \& Toumi, N.: Finite elements in
Archimedean $\Phi$-algebras. Positivity 29, 55 (2025).

\bibitem {10}Luxemburg, W.A.J., Zaanen, A.C.: Riesz Spaces I. North-Holland,
Amsterdam (1971).

\bibitem {11}Schaefer, H.H.: Banach Lattices and Positive Operators. Springer,
Berlin, Heidelberg, New York (1974).

\bibitem {12}Toumi, M. A.: Characterization of self-majorizing elements in
Archimedean unital $f$-algebras. Arch. Math. 123, 253--261 (2024).

\bibitem {13}Toumi, M.A., Jbeli, Z.: Characterizations of finite and
self-majorizing functions in $C(X)$. Anal Math (2025).

\bibitem {14}Toumi, M.A.: Orthogonally additive polynomials on Dedekind
$\sigma$-complete vector lattices, Math. Proc. R. Ir. Acad. 110A, no. 1, 83-94 (2010).

\bibitem {15}Toumi, M.A.: Self-majorizing elements in Banach lattice unital
$f$-algebras. Positivity 28, 41 (2024).

\bibitem {16}Toumi, M.A., Toumi, N. \& Jbeli, Z.: Centered orthomorphisms and
the Wickstead problem. Positivity 30, 29 (2026).

\bibitem {17}Weber, M. R.: Finite Elements in Vector Lattices. Berlin-Boston,
W. de Gruyter GmbH, (2014).

\bibitem {18}Wickstead, A.W.: Representation and duality of multiplication
operators on Archimedean Riesz spaces, Compos. Math. 35(3), 225--238 (1977).
\end{thebibliography}
\end{document}